\documentclass [11pt] {article}

\usepackage{amsmath,amssymb}

\newtheorem{theorem}{Theorem}

  \newtheorem{corollary}{Corollary}

  \newcommand {\qed} {\null \hfill \rule{2mm}{2mm}}

\def\sep{{\rm sep}}

\begin {document}

\title{{\Large{\bf Root separation for irreducible integer polynomials}}}

\author{Yann Bugeaud and Andrej Dujella}

\date{}
\maketitle

\footnotetext{
 \vspace{1ex} \\
The authors were supported by the French-Croatian bilateral COGITO
project {\it Diophantine approximations}.}

\begin{abstract}
We establish new results on root separation of integer, irreducible
polynomials of degree at least four. These improve earlier bounds of
Bugeaud and Mignotte (for even degree) and of Beresnevich, Bernik, and
G\"otze (for odd degree).
\end{abstract}

\section{Introduction}

The height $H(P)$ of  an integer
polynomial $P(x)$ is
the maximum of the absolute values of its
coefficients.
For a separable
integer polynomial $P(x)$ of degree $d \ge 2$
and with distinct roots $\alpha_1, \ldots , \alpha_d$,
we set
$$
\sep(P) = \min_{1 \le i < j \le d} \, |\alpha_i - \alpha_j|
$$
and define $e(P)$ by
$$ \sep(P) = H(P)^{-e(P)}. $$
Following the notation from \cite{BM2}, for $d \ge 2$, we set
$$
e(d) := \limsup_{{\rm deg}(P) = d, H(P) \to + \infty} e(P)
$$
and
$$
e_{{\rm irr}} (d) := \limsup_{{\rm deg}(P) = d, H(P) \to + \infty} e(P),
$$
where the latter limsup is taken over the irreducible integer polynomials
$P(x)$ of degree $d$.
A classical result of Mahler \cite{Mah64}
asserts that $e(d) \le d-1$ for all $d$, and
it is easy to check that $e_{{\rm irr}}(2)=e(2)=1$. There is only one other value
of $d$ for which $e(d)$ or $e_{{\rm irr}}(d)$ is known, namely $d=3$, and
we have $e_{{\rm irr}}(3)=e(3)=2$, as proved, independently, by Evertse \cite{Ev}
and Sch\"onhage  \cite{Sch}.
For larger values of  $d$, the following lower bounds have been established
by Bugeaud and Mignotte in \cite{BM1}:
$$ e_{{\rm irr}}(d) \geq d/2, \quad \mbox{for even $d\geq 4$}, $$
$$ e (d)  \geq (d+1)/2, \quad \mbox{for odd $d\geq 5$}, $$
$$ e_{{\rm irr}}(d)  \geq (d+2)/4, \quad \mbox{for odd $d\geq 5$}, $$
while Beresnevich, Bernik, and G\"otze \cite{BBG} proved that
$$ e_{{\rm irr}}(d)\geq (d+1)/3, \quad \mbox{for every $d\geq 2$}. $$
Except those from \cite{BBG}, the
above results are obtained by presenting explicit families
of (irreducible) polynomials of degree $d$ whose roots are close enough.
The ingenious proof in \cite{BBG} does not give any explicit example of such polynomials, but
shows that algebraic numbers of degree $d$ with a close conjugate
form a `highly dense' subset in the real line.

The aim of the present note is to improve all known lower bounds for $e_{{\rm irr}}(d)$
when $d\geq 4$.

\begin{theorem} \label{tm1}
For any integer $d\geq 4$, we have
$$ e_{{\rm irr}}(d)  \geq \frac{d}{2} + \frac{d-2}{4(d-1)}. $$
\end{theorem}

To prove Theorem 1, we construct explicitly,
for any given degree $d \ge 4$, a one-parametric family of irreducible integer
polynomials $P_{d,a}(x)$  of degree $d$. We postpone to Section 3 our general
construction and give below some numerical examples in small degree.

For $a \ge 1$, the roots of the polynomial
$$ P_{4,a}(x)=(20a^4-2)x^4 +(16a^5+4a)x^3 +(16a^6+4a^2)x^2 +8a^3x+1, $$
are approximately equal to:
\begin{eqnarray*}
 r_1 &=& -(1/4) a^{-3} - (1/32) a^{-7} - (1/256) a^{-13} + \ldots , \\
 r_2 &=& -(1/4) a^{-3} - (1/32) a^{-7} + (1/256) a^{-13} + \ldots , \\
 r_3 &=& -(2/5) a + (11/100) a^{-3} + (69/4000) a^{-7} + (4/5) a \,i + \ldots ,\\
 r_4 &=& -(2/5) a + (11/100) a^{-3} + (69/4000) a^{-7} - (4/5) a \,i + \ldots .
\end{eqnarray*}
Since $H(P_{4,a})=O(a^6)$ and $\sep(P_{4,a})=|r_1-r_2|=O(a^{-13})$,
we obtain by letting $a$ tend to infinity that $e_{{\rm irr}}(4) \geq 13/6$.

A similar construction for degree five gives the family of polynomials
$$ P_{5,a} (x) = (56a^5-2)x^5+(56a^6+4a)x^4+(80a^7+4a^2)x^3+
(100a^8+8a^3)x^2+20a^4x+1 $$
with two close roots
$$
(1/10) a^{-4}+(1/250) a^{-9}+(3/25000) a^{-14}-(3/250000) a^{-19}
\pm (\sqrt{10}/500000) a^{-43/2} + \ldots ,
$$
and we obtain that $e_{{\rm irr}}(5)\geq 43/16$.

Our construction is applicable as well for $d=3$.
It gives the family
$$ P_{3,a} (x) =(8a^3-2)x^3+(4a^4+4a)x^2+4a^2x+1 $$
with close roots
$-(1/2) a^{-2}-(1/4) a^{-5} \pm (\sqrt{2}/8) a^{-13/2}$,
showing that $e_{{\rm irr}}(3) \geq 13/8$.
%By the way, this degree 3 example gives exponent $0.75$ in Hall's conjecture:
%$X=4n^4-12n$, $Y=8n^6-36n^3+27$, $|X^3-Y^2|=216n^3-729$.
This is weaker than the known result $e_{{\rm irr}}(3)=2$,
but it could be noted that in the examples showing that
$e_{{\rm irr}}(3)=2$ the coefficients of
the polynomials involved have exponential growth, while
in our example the coefficients have polynomial growth, only.

The constant term of every polynomial $P_{d,a}(x)$
constructed in Section 3 is equal to $1$.
This means that the reciprocal polynomial of $P_{d,a}(x)$ is monic.
Therefore, Theorem \ref{tm1} gives also a
lower bound for the quantity
$$
e_{{\rm irr}}^*(d) := \limsup_{{\rm deg}(P) = d, H(P) \to + \infty} e(P),
$$
where the limsup is taken over the {\it monic} irreducible integer polynomials.
Regarding this quantity,
the following estimates have been established
by Bugeaud and Mignotte in \cite{BM2}:
$$ e_{{\rm irr}}^*(2) = 0, \quad e_{{\rm irr}}^*(3) \geq 3/2, $$
$$ e_{{\rm irr}}^*(d) \geq (d-1)/2, \quad \mbox{for even $d\geq 4$}, $$
$$ e_{{\rm irr}}^*(d)  \geq (d+2)/4, \quad \mbox{for odd $d\geq 5$}, $$
while Beresnevich, Bernik, and G\"otze \cite{BBG} proved that
$$ e_{{\rm irr}}^*(d)\geq d/3, \quad \mbox{for every $d\geq 3$}. $$
In particular, for $d=5$, the current best estimate is $e_{{\rm irr}}^*(5) \ge 7/4$.

Our construction allows us to improve these results when $d$ is odd and
at least equal to $7$.

\begin{theorem} \label{tm2}
For any odd integer $d\geq 7$,  we have
$$
e_{{\rm irr}}^*(d) \geq \frac{d}{2} + \frac{d-2}{4(d-1)} -1.
$$
\end{theorem}

To prove Theorem 2, we simply observe that if $\alpha$ and $\beta$ denote the
two very close roots of a polynomial $P_{d,a} (x)$ constructed
in Section 3, then $\alpha$ and $\beta$ satisfy
$$
|\alpha |^{-1}, |\beta |^{-1} = O(a^{d -1}) = O (H(P_{d, a})^{1/2}),
$$
and
$$
\left| \frac{1}{\alpha} -  \frac{1}{\beta} \right|
= \frac{|\alpha - \beta |}{\alpha \beta}
$$
is very small, where, clearly, $1/\alpha$ and $1/ \beta$ are roots
of the reciprocal polynomial of $P_{d,a} (x)$.

\section{Application to Mahler's and Koksma's classifications of numbers}

The families of polynomials constructed for the proof of Theorem 1
can be used in the context of \cite{Bu1}.
Let $d$ be a positive integer. Mahler and, later, Koksma, introduced
the functions $w_d$ and $w_d^*$
in order to measure the quality of approximation
of real numbers by
algebraic numbers of degree at most $d$.
For a real number $\xi$, we denote by
$w_d (\xi)$ the supremum of the exponents $w$ for which
$$
0 < |P(\xi)| < H(P)^{-w}
$$
has infinitely many solutions in integer polynomials $P(x)$ of
degree at most $d$.
Following Koksma, we denote by
$w_d^* (\xi)$ the supremum of the exponents $w^*$ for which
$$
0 < |\xi - \alpha| < H(\alpha)^{-w^*-1}
$$
has infinitely many solutions in real algebraic numbers $\alpha$ of
degree at most $d$. Here, $H(\alpha)$ stands for the na\"\i ve height of $\alpha$,
that is, the na\"\i ve height of its minimal defining polynomial.

For an overview of results on $w_d$ and $w_d^*$,
the reader can consult \cite{BuLiv}, especially Chapter 3.
Let us just mention
that it is quite easy to establish the inequalities
$$
w_d^*(\xi) \le w_d (\xi) \le w_d^*(\xi) + d -1,
$$
for any transcendental real number $\xi$,
and that
$$
w_d^*(\xi) = w_d (\xi) = d
$$
holds for almost all real numbers $\xi$, with respect to the
Lebesgue measure.

For $d \ge 2$, R. C. Baker \cite{Bak} showed
that the range of values of the function $w_d - w_d^*$ includes the
interval $[0, (d-1)/d]$. This has been substantially improved in
\cite{Bu1}, where it is shown that
the function $w_d - w_d^*$ can take any value
in $[0, d/4]$. Further results are obtained in \cite{Bu2,Bu3},
including that the function $w_2 - w_2^*$ (resp. $w_3 - w_3^*$)
takes any value in $[0, 1)$ (resp. in $[0, 2)$).
The proofs in \cite{Bu1,Bu2,Bu3} make use of families of polynomials
with close roots. In particular, the upper bound $d/4$ is obtained
by means of the family of polynomials $x^d - 2 (a x - 1)^2$ of height
$2 a^2$ and having two roots separated by $O(a^{-(d+2)/2})$.

\begin{corollary} \label{cor1}
For any integer $d \geq 4$,  the function $w_d - w_d^*$
takes every value in the interval
$$
\left[0, \frac{d}{2} + \frac{d-2}{4(d-1)}\right).
$$
\end{corollary}

This corollary is established following the main lines of the proofs
of similar results established in \cite{Bu1,Bu2,Bu3}.
We omit the details.

\section{Proof of Theorem 1: construction of families of integer polynomials}

For each integer $d\geq 3$, we construct a sequence of integer polynomials $P_{d,a}(x)$
of degree $d$ and arbitrarily large height  having two roots very close to each other,
and whose coefficients are polynomials in the
parameter $a$.

For $i \ge 0$, let $c_i$ denote the $i$th Catalan number defined by
$$
c_i=\frac{1}{i+1}{2i \choose i}.
$$
The sequence of Catalan numbers $(c_i)_{i\geq 0}$ begins as
$$
1,1,2,5,14,42,132,429,1430,\ldots
$$
and satisfies the recurrence relation
\begin{equation} \label{cat}
c_{i+1}=\sum_{k=0}^{i} c_k c_{i-k},
\quad \hbox{for $i \ge 0$}.
\end{equation}
For integers $d \ge 3$ and $a \ge 1$, consider the polynomial
\begin{eqnarray*}
  P_{d,a}(x) &=& (2c_0 ax^{d-1}+2c_1 a^2 x^{d-2}+ \ldots +2c_{d-2} a^{d-1} x)^2  \\
  & & \mbox{} - (4c_1 a^2x^{2d-2} +4c_2 a^3x^{2d-3}+ \ldots +4c_{d-2} a^{d-1} x^{d+1}) \\
&  & \mbox{} + (4c_1 a^2 x^{d-2}+4c_2 a^3 x^{d-3}+ \ldots +4c_{d-2} a^{d-1} x ) \\
& & \mbox{}+ 4a x^{d-1} - 2x^d +1,
\end{eqnarray*}
which generalizes the polynomials $P_{3,a}(x)$, $P_{4,a}(x)$, $P_{5,a}(x)$
given in Section 1.
It follows from the recurrence (\ref{cat}) that $P_{d,a}(x)$ has degree
exactly $d$, and not $2d-2$,
as it seems at a first look.
Furthermore, we check that the height of $P_{d,a}(x)$
is given by the coefficient of $x^2$, that is,
$$
H(P_{d,a}) = 4 c^2_{d-2} a^{2d-2} + 4 c_{d-3} a^{d-2}.
$$

By applying the Eisenstein criterion with the prime $2$ on the reciprocal polynomial
$x^d P_{d,a}(1/x)$, we see that the polynomial $P_{d,a}(x)$ is irreducible. Indeed,
all the coefficients of $P_{d,a}(x)$ except the constant term are even,
but its leading coefficient,
which is equal to $4c_{d-1}a^d -2$, is not divisible by 4.

Writing
$$
g=g(a,x) = 2 c_0 a x^{d-1} + 2 c_1 a^2 x^{d-2} + \ldots + 2 c_{d-2} a^{d-1} x,
$$
we see that
$$
P_{d,a} (x) = (1 + g)^2 + x^d (4 a x^{d-1} - 2 (1 + g)).
$$
Rouch\'e's theorem shows that $P_{d,a}(x)$ has exactly two roots
in the disk centered at the origin and of radius $1/2$.
Clearly, $(1 + g)^2$ has a double root, say $x_0$,
close to $- 1/(2  c_ {d-2} a^{d-1})$. More precisely,
we have
$$
x_0 = -  a^{-d+1} /(2  c_ {d-2}) + O (a^{-2d+1}).
$$
Here and below, the numerical constants implied in $O$
are independent of $a$.

The polynomial $P_{d,a}(x)$ has two distinct roots
close to $x_0$, since
the term $x^d (4 a x^{d-1} - 2 (1 + g))$ is a small perturbation when $x$ is near $x_0$.
Below we make this more precise.

Observe that, for any real number $\delta$, we have
$$
1 + g \bigl(a, x_0 + \delta a^{-d^2+d/2+1} \bigr)  = 2 \delta c_{d-2}
a^{-d^2+ 3d/2} + O(a^{-d^2+d/2} ),
$$
thus
$$
P_{d,a} (x_0 + \delta a^{-d^2 + d/2 +1})
= 4 \bigl(\delta^2 c_{d-2}^2  - (2 c_{d-2})^{-2d+1} \bigr)
a^{-2d^2+3d} + O(a^{- 2d^2 + 5d/2}).
$$
Let
$$
\delta_0=\frac{1}{2^{d-1/2} c_{d-2}^{d+1/2}}.
$$
Then for every sufficiently small
$\varepsilon >0$
and every sufficiently large $a$, we have
$$
P_{d,a}(x_0 \pm (\delta_0+\varepsilon) a^{-d^2+d/2+1}) > 0
$$
and
$$
P_{d,a}(x_0 \pm (\delta_0-\varepsilon) a^{-d^2+d/2+1}) < 0.
$$
This shows that $P_{d,a}(x)$ has a root $x_1$ in the interval
$$
\bigl(x_0 - (\delta_0+\varepsilon) a^{-d^2+d/2+1}, x_0 -
(\delta_0-\varepsilon) a^{-d^2+d/2+1} \bigr)
$$
and a root $x_2$
in the interval
$$
\bigl(x_0 + (\delta_0-\varepsilon) a^{-d^2+d/2+1}, x_0 +
(\delta_0+\varepsilon) a^{-d^2+d/2+1} \bigr).
$$
This yields
$$
{\rm sep}(P_{d,a}) \leq 2 (\delta_0 + \varepsilon) a^{-d^2+d/2+1}.
%=\frac{1}{2^{d-3/2} c_{d-2}^{d+1/2}a^{d^2-d/2-1}}.
$$
Since $H(P_{d,a})=O(a^{2d-2})$, this gives
$$
e_{{\rm irr}}(d)\geq \frac{2d^2-d-2}{4(d-1)} = \frac{d}{2} + \frac{d-2}{4(d-1)},
$$
by fixing the arbitrarily small positive
real number $\varepsilon$ and letting $a$ tend to infinity.
The proof of Theorem \ref{tm1} is complete.
\qed

\bigskip

{\small \noindent
Yann Bugeaud \\
Universit\'e de Strasbourg  \\
D\'epartement de Math\'ematiques \\ 7, rue Ren\'e Descartes \\
67084 Strasbourg, France \\
{\em E-mail address}: {\tt bugeaud@math.unistra.fr}}

\bigskip

{\small \noindent
Andrej Dujella \\
Department of Mathematics \\ University of
Zagreb
\\ Bijeni\v{c}ka cesta 30 \\
10000 Zagreb, Croatia \\
{\em E-mail address}: {\tt duje@math.hr}}


\medskip




\bigskip

\begin{thebibliography}{99}

\bibitem{Bak}
R. C. Baker,
{\it On approximation with algebraic numbers of bounded degree},
Mathematika 23 (1976), 18--31.



\bibitem{BBG}
V. Beresnevich, V. Bernik, and F. G\"otze,
{\it The distribution of close conjugate algebraic numbers},
Compositio Math. 146  (2010), 1165--1179.


\bibitem{Bu1}
Y. Bugeaud,
{\it Mahler's classification of numbers compared with Koksma's},
Acta Arith. 110 (2003), 89--105.


\bibitem{Bu2}
Y. Bugeaud,
{\it Mahler's classification of numbers compared with Koksma's, III},
Publ. Math. Debrecen {65} (2004), 305--316.


\bibitem{BuLiv}
Y. Bugeaud,
Approximation by algebraic numbers.
Cambridge Tracts in Mathematics, Cambridge, 2004.

\bibitem{Bu3}
Y. Bugeaud,
{\it Mahler's classification of numbers compared with Koksma's, II}.
In: Diophantine approximation,  107--121,
Dev. Math., 16, Springer Wien New York, Vienna, 2008.

\bibitem{BM1}
Y. Bugeaud and M. Mignotte,
{\it On the distance between roots of integer polynomials},
Proc. Edinburgh Math. Soc. 47 (2004), 553--556.


\bibitem{BM2}
Y. Bugeaud and M. Mignotte,
{\it Polynomial root separation},
Intern. J. Number Theory 6 (2010), 587--602.

\bibitem{Ev}
J.-H. Evertse,
{\it Distances between the conjugates of an algebraic number},
Publ. Math. Debrecen 65 (2004), 323--340.


\bibitem{Mah64}
K. Mahler,
{\it An inequality for the discriminant of a polynomial},
Michigan Math. J. 11 (1964), 257--262.

\bibitem{Sch}
A. Sch\"onhage,
{\it Polynomial root separation examples},
J. Symbolic Comput.  41  (2006),  1080--1090.



\end{thebibliography}
\end{document}